\documentclass[12pt,a4paper]{article}
\usepackage[utf8]{inputenc}
\usepackage[T1]{fontenc}
\usepackage{graphicx}
\usepackage{grffile}
\usepackage{longtable}
\usepackage{wrapfig}
\usepackage{rotating}
\usepackage[normalem]{ulem}
\usepackage{amsmath}
\usepackage{textcomp}
\usepackage{amssymb}
\usepackage{capt-of}
\usepackage[colorlinks=true]{hyperref}
\usepackage{minted}
\usepackage[margin=1.5cm,includefoot]{geometry}
\usepackage{amsmath,amssymb}

\author{Siniša Miličić}
\date{21 November 2026}
\title{Questions, Explanation and Synthesis in Service Mathematics\\\medskip
\large Course Design with Didactic Units and Hybrid Teaching}
\hypersetup{
 pdfauthor={Siniša Miličić},
 pdftitle={Questions, Explanation and Synthesis in Service Mathematics},
 pdfkeywords={service mathematics; student explanation; teacher synthesis; formative feedback; didactic units; mathematical compression},
 pdfsubject={},
 pdfcreator={Emacs 30.2 (Org mode 9.7.30)}, 
 pdflang={English}}
\usepackage[style=authoryear,backend=biber]{biblatex}
\begin{document}

\maketitle
\begin{abstract}
Students may have access to a convincing mathematical explanation without being able to reconstruct or adapt it. A course must provide opportunities to develop that understanding while keeping the work sustainable for students and teachers. This article develops a methodology for service mathematics in which a question ledger guides both course authoring and classroom decisions. Course-relative didactic units organise explanatory purposes, prerequisites and connections, initially calibrated by a five-to-seven-minute expert account. Teaching combines preparation, student exposition and supported repair with substantial teacher exposition and synthesis. Their allocation changes with the topic, audience and conditions of teaching. Compression and recovery explain how developed constructions become reusable while their grounds remain accessible. Three analyses from mathematics for computing—full Sudoku encoding, grayscale projection and a software failure contract—show how contrasting cases guide explanation and revision. The practice-informed contribution connects coherent course design with adaptable teaching and manageable preparation, and identifies questions for evaluating learning, participation and sustained teacher effort.
\end{abstract}
\section{Introduction}
\label{sec:org32ce310}

Human understanding is the main concern of good mathematics teaching. Students need to formulate situations mathematically, explain why a construction works and recognise what changes when its assumptions change. A course must make this work possible within the time available. The practical questions are connected: why is this content here, what should students be able to explain, and what conditions, resources and activities enable them to do that work?

The aim is \emph{human mathematical fluency}: movement between particular cases, representations, compact objects and general arguments. In service mathematics, questions from another discipline provide multiple starting points. Mathematics then has to become a coherent account that students can use and reopen. A list of applications attached to an otherwise unchanged syllabus does not by itself accomplish this.

Developed through the author's teaching of service mathematics and computing, the methodology coordinates course authoring with opportunities for students to prepare, explain and repair mathematical accounts. Teacher exposition and synthesis establish standards and connect the work. Didactic units specify explanatory purposes and prerequisites; workload records connect them to feasible activities.

The central argument is that explicit explanatory obligations can coordinate decisions across course, unit and session scales. The worked analyses examine this coordination through changes in assumptions, representations and dependencies; evaluation of educational effects remains future work.

Ready access to explanations, including those generated by large language models, makes the distinction between possessing and understanding an account especially visible. Live dialogue provides opportunities to reconstruct an argument, interpret a symbol and respond to a changed case. This gives dialogue both a formative and an evidential role.

Table \ref{tab:org1e2012e} locates the method's three planning scales. The question ledger connects them: classroom difficulties can prompt changes to prerequisites, unit boundaries or course questions. These are interdependent decisions, with feedback from teaching to authoring.

\begin{table}[tbp]
\caption{\label{tab:org1e2012e}Course, unit and session planning, with feedback between scales.}
\centering
\begin{tabular}{p{0.12\linewidth}p{0.38\linewidth}p{0.35\linewidth}}
Scale & Design decisions & Evidence for revision\\
\hline
Course & Realism lenses, organising questions, outcomes and scope boundaries. & Unjustified content, missing connections or unclear scope boundaries.\\
Unit & Explanatory purpose, prerequisites, contrasting cases and 5–7-minute expert calibration. & Hidden dependencies, rushed accounts or several independent purposes.\\
Session & Advance preparation; group preparation where appropriate; exposition and repair; protected teacher synthesis. & Explanations, peer contributions, recurring difficulties and time used.\\
\end{tabular}
\end{table}
\section{Preparation, explanation and teacher synthesis}
\label{sec:org905f7a1}

\subsection{A hybrid teaching arrangement}
\label{sec:orgcdb18b0}

The central student task is an \emph{expositional problem}: prepare an account that another person can question. The task may concern a theorem, construction, example, experiment, motivation or piece of code. Materials, prerequisites, expected depth and organising questions are available in advance. Texts, diagrams, examples, code and recordings can support preparation according to the task.

In the author's usual full-time configuration, teaching takes place in a double forty-five-minute block. Students review material beforehand and prepare in groups in class. Students from the groups then present, with peers and the teacher participating in the discussion; about five presentations usually fit the block. Online part-time students may prepare individually before a live exchange.

Following others' presentations is part of every student's work. Peers can check a claim against their preparation, identify an omitted case, request a justification or connect accounts. Their additions and discussion can contribute to assessment against the relevant outcomes. Mathematical substance matters; frequency or speed of speaking alone does not establish attainment.

Partially open presentation tasks give students room to choose and develop a topic under a stated thematic or depth constraint. They are particularly useful in the online part-time configuration. Student creativity can introduce a lively application or a tangential connection that enriches the course. Such presentations may run longer and encompass several units; their assessment remains tied to the relevant outcomes.

The arrangement is \emph{hybrid} because it combines student exposition, dialogue and direct teacher exposition. It is \emph{parametric} because their allocation, together with preparation, grouping, depth and support, can be adjusted to the course. Otten and colleagues' account of structural variation in flipped mathematics instruction provides a precedent for specifying these choices \autocite{OttenEtAl2023Variation}.

The allocation also varies by topic. In the author's introductory teaching for informatics and computer science students, the full Sudoku development is taught frontally, with all examples supplied. A coherent expert account preserves the recognisable problem while making its formal construction accessible. Kurnik's \emph{problemska} teaching, organised around student investigation, and \emph{heuristička} teaching, supported by teacher questioning, offer related approaches to adjusting guidance \autocite{Kurnik2002Problemska,Kurnik2006Heuristicka}. Direct exposition likewise has a deliberate place.

Self-explanation and learning-by-teaching research motivate constructing and communicating accounts \autocite{ChiEtAl1994SelfExplanations,FiorellaMayer2013Teaching}. Retrieval, inference and response to questions can contribute differently; experiments identify conditions in which retrieval explains a learning-by-teaching benefit \autocite{KohEtAl2018Retrieval}. Here students reconstruct relationships and make them available for discussion, using notes and representations where useful.
\subsection{Feedback, repair and assessment}
\label{sec:org8854056}

Questions during a presentation should locate the next useful educational action. Asking what a variable denotes, requesting an example and asking why an argument covers every case address different difficulties. A contrasting case can expose a missing condition; an unfamiliar representation may first require direct explanation. The teacher may supply a prerequisite or invite a peer contribution, then return the explanation to the presenter. The student must articulate the repaired connections and justifications: a teacher's completed solution cannot substitute for this work.

Feedback matters when it changes the work. The student revises an explanation, peers reconsider a construction, and the teacher adjusts a prompt or the sequence. Repeated difficulties can later motivate changes to materials, prerequisites or unit boundaries. This connects immediate dialogue with the formative assessment functions of supporting learning and informing teaching \autocite{NicolMacfarlaneDick2006Feedback,BlackWiliam2009}. Feedback is not uniformly beneficial: its focus and use matter, as the heterogeneous effects reviewed by Kluger and DeNisi show \autocite{KlugerDeNisi1996Feedback}. Correction should therefore identify the mathematical issue and lead to an opportunity to act on it.

Presentations are graded on the explanation achieved \emph{after supported repair}, against stated course outcomes and conditions of support. Two students reaching the same required level warrant the same judgement of that attainment, even if they needed different assistance. Prompt counts are not a hidden grading criterion. The student must nevertheless articulate the repaired account. Dynamic assessment provides a related approach to integrating assistance and assessment \autocite{LantolfPoehner2011Dynamic}.

Higher attainment can involve deeper justification, adaptation, integration or a creative project within the course's scope. These opportunities develop students' capacities beyond a required threshold. Common criteria should distinguish such achievement from the amount of assistance received. The supported presentation contributes evidence for its specified outcomes; any exams or integrative projects assess their own stated outcomes and conditions. Presentation style counts where it expresses an outcome, such as communicating a software pattern through code.
\subsection{What teacher synthesis contributes}
\label{sec:org2da814c}

Teacher synthesis needs protected time. It establishes which assumptions, notation and justifications should govern the account, resolves consequential errors and connects contributions. Students should hear an expert demonstrate how much detail is needed, which distinctions matter and how an explanation can become concise without losing its grounds.

Synthesis also prepares mathematical compression. Several cases can be organised by one relation; a developed process can become an object available for further reasoning. The teacher identifies what is being retained, where it applies and how its interior can be recovered. Connections and syntheses can themselves be planned as didactic units or combinations of units.

Synthesis may develop throughout a teacher-led exposition or follow student contributions. The teacher can prepare its essential connections and adapt them to the difficulties encountered. Its scope and timing remain planned. Necessary exposition can require fewer presentations or a predominantly teacher-led session; protecting synthesis constrains the other activities.

Brousseau's account of \emph{institutionalisation} supplies a close precedent: the teacher gives selected classroom knowledge a shared mathematical status and organises what students should retain \autocite{Brousseau2006Research}. That public organisation and an individual's encapsulation are related but distinct. A coherent synthesis makes a reusable structure available; hearing it does not establish that every student can already use or reconstruct it.

Teacher exposition has established places within flipped mathematics teaching. Feudel and Fehlinger's lecture-oriented flipped analysis course uses advance reading and classroom expert modelling of definitions, proofs and connections \autocite{FeudelFehlinger2023LectureOriented}. Maciejewski combines preparation, peer discussion and teacher resolution in calculus \autocite{Maciejewski2016Flipping}. The present proposal connects those functions to student explanation, repair, question-led authoring and unit-based resource planning.
\section{Questions that justify and organise the course}
\label{sec:org6cdc707}

\subsection{Three realism lenses}
\label{sec:orgf92a31e}

An organising question should explain why content belongs in the course. Compatible schedules and plans justify propositional logic in introductory computing mathematics. How to solve a Sudoku by programming offers an experiential entry into related constructions. The course must establish what transfers and which assumptions differ.

Three realism lenses help construct and examine such questions. \emph{Professional realism} concerns purposes, tasks and constraints in actual or anticipated work. \emph{Disciplinary realism} concerns how a field teaches and communicates among peers: its representations, explanatory conventions and ways of posing questions. These practices may be mathematical to a degree even when practitioners do not describe them as mathematics. \emph{Experiential realism} concerns what students can imagine, interpret and act upon from their experience.

The lenses are structural as well as analytical: they help construct and examine tasks. The teacher checks a question's purpose, its disciplinary account and this audience's access to it. A professionally plausible question can be unfamiliar to students; a familiar puzzle can need a further connection to professional work. These overlapping lenses locate what the course must develop without constituting an exhaustive taxonomy.

One task can meet all three concerns as students become familiar with its purpose and representations. Disciplinary communication resembles the differentiated practices studied in disciplinary literacy \autocite{ShanahanShanahan2008Disciplinary}, while experiential realism can include familiar mathematical objects. Throughout this development, professional questions are worked on as mathematics: their assumptions, constructions and conclusions require mathematical explanation and justification.

Realistic Mathematics Education (RME) explains how meaningful activity can develop into mathematical representations and more general reasoning \autocite{Freudenthal1991Revisiting,GravemeijerDoorman1999}. The Anthropological Theory of the Didactic (ATD) helps organise the relation between questions, techniques and their justifications within institutions \autocite{Chevallard2019Introducing}. These perspectives prevent the organising question from becoming a decorative introduction. It should influence the mathematical work, the explanations required and the depth pursued. Existing syntheses of inquiry, RME and ATD establish a substantial background for this orientation \autocite{ArtigueBlomhoj2013Inquiry,ArtigueEtAl2020Trajectories}.
\subsection{The question ledger and directions of authoring}
\label{sec:orgc0664b9}

Test-driven development (TDD) inspires a course-authoring procedure: formulate a question the account does not yet answer, construct an adequate answer, then refactor the account while preserving its obligations \autocite{Beck2003TDD}. Refactoring may align notation, remove duplication, reorder material or introduce an intermediate construction. Rechecking the affected questions exposes unsupported transitions and misleading examples.

The questions test course purposes, integration and local constructions. Table \ref{tab:org423292a} presents these levels from course purposes downwards. Authoring can begin with a syllabus, grow from a worked construction or proceed between these levels. The ledger tests the authored account: some mathematical checks are executable, while explanatory adequacy requires judgement.

\begin{table}[tbp]
\caption{\label{tab:org423292a}Proposed question-ledger excerpt for the teacher-led propositional-logic development and its wider course purposes.}
\centering
\begin{tabular}{p{0.10\linewidth}p{0.28\linewidth}p{0.28\linewidth}p{0.20\linewidth}}
Level & Question & Contrasting case or test & Target unit or connection\\
\hline
Course & Does the account support constructing and justifying a model for another compatible-choice problem? & A workflow allowing idling; identify changed assumptions and supplied dependencies. & Transfer, course synthesis and scope.\\
Integration & Why do the combined restrictions describe precisely the valid Sudoku completions? & A filling valid by rows and columns but invalid by boxes; justify both directions. & Local constraints and global correctness.\\
Integration & Why does the derived CNF preserve the intended meaning? & An assignment separating an incorrect conversion from the original formula. & Formula structure and equivalence transformations.\\
Local & Which restriction requires at least one, at most one or exactly one choice? & Zero, one and two selected choices. & Cardinality constructions.\\
Local & What do the variables and their truth values say about a configuration? & Hold the formula fixed and change its assignment. & Propositions, formulas and interpretations.\\
\end{tabular}
\end{table}

The \emph{question ledger} serves both advance authoring and live teaching. It records what the account must explain and helps the teacher locate the current question, revisit an assumption or defer an extension. Course purposes can appear at syllabus level; the teacher selects questions, contrasts and expectations to share for student preparation. A compact view of the same record can support a session without duplicating the full course plan.

Classroom difficulties inform subsequent revision. For example, transferring exactly-one constraints to scheduling fails the ledger's idling test when inactivity is permitted. Refactoring separates required selection from mutual exclusion, then checks both applications again.

This procedure is related to constructive alignment and earlier assessment-driven design \autocite{Biggs1996Alignment,SchwiegerSurendran2013AssessmentDriven}. The questions themselves must be examined for omissions, inappropriate scope and weak connections. Refactoring should produce a coherent mathematical story, with connections and syntheses represented by units whose purposes and dependencies are explicit.

\emph{Middle-out design} is one route: start from a workable construction and choose how far to open its dependencies and extensions. In an introductory logic course, students may explain the correspondence between a formula and a planning problem while using a supplied solver for Boolean satisfiability (SAT). A senior course can investigate its algorithms, complexity, implementation and testing. ATD's treatment of black and clear boxes provides a related account of selective investigation \autocite{BoschEtAl2020Casebook}. Outcomes, prerequisites and mathematical coherence determine the boundary: not every available direction needs investigation in the present course.
\section{Mathematical objects, compression and didactic units}
\label{sec:orgf54a5a1}

\subsection{Objects that can be used and reopened}
\label{sec:org274282f}

Davis and Hersh characterise mathematics through mental objects with reproducible properties \autocite{DavisHersh1981MathematicalExperience}. This is a useful orientation for teaching: students need objects whose behaviour they can examine, communicate and reconstruct. Borovik develops the importance of compression and recovery in mathematical practice \autocite{Borovik2010Microscope}. A compact object can carry considerable prior work, provided that suitable representations and arguments make that work accessible again.

Freudenthal approaches a related issue didactically. Meaningful use and the constitution of mental objects can precede explicit concept attainment; mathematical structures organise phenomena encountered in activity \autocite{Freudenthal1983Phenomenology}. Davis and Hersh offer a characterisation of mathematics, Borovik examines mathematical practice, and Freudenthal asks how objects become meaningful for learners. Their combination directs attention both to reusable structures and to the situations from which students can develop them.

\emph{Encapsulation} is central here: a process developed through several steps becomes available as an object in subsequent thought. Process–object accounts, including Sfard's discussion of reification, make this transition explicit \autocite{Sfard1991DualNature,GrayTall1994}. Borovik also cautions against reducing mathematical compression to reification alone. New representations, languages and ways of organising arguments can compress work without fitting one process-to-object transition.

The exactly-one construction illustrates why compact notation, acquired understanding and curricular scope need separate attention. Naming the construction abbreviates a formula. A learner develops a reusable component by constructing and justifying that formula, then adapting it when the requirement changes—for example, when selecting no option becomes permissible. The same learner may use a supplied SAT solver without studying its search algorithm. The first issue concerns representation, the second what the learner can reconstruct, and the third what the course deliberately leaves unopened. A short expression alone does not tell the teacher which work remains to be learned.
\subsection{Directions of mathematical thought}
\label{sec:org9b7b99e}

Kurnik's accounts of scientific and heuristic methods provide a planning vocabulary: analysis and synthesis, abstraction and concretisation, generalisation and specialisation, induction and deduction, and analogy \autocite{Kurnik2008Znanstvenost,Kurnik2006Heuristicka}. Within a longer heuristic tradition \autocite{Polya1954Plausible}, these directions suggest ways to develop and reopen mathematical objects. The \emph{compass} is the present article's metaphor for their use in planning.

A further emphasis is \emph{cataloguing}: identifying, gathering and differentiating particulars. It supplies cases and distinctions for the ledger. Analogy connects through commonalities and structural correspondences, helping synthesis identify reusable constructions. Their opposed emphases interact: recognised common structure can make a difference visible \autocite{GentnerGunn2001Differences}. The examples below use this distinction to examine choice configurations, transformations and failure traces.

Encapsulation, analogy, induction, abstraction and generalisation often compress mathematical work, although an inductively suggested pattern still needs justification and an analogy has limits. Cataloguing, specialisation and concretisation help inspect what an account preserves or overlooks. The compass identifies useful instructional moves; it neither ranks attainment nor requires every example to traverse every direction.
\subsection{Units and the cost of learning}
\label{sec:orge82fca8}

A \emph{didactic unit} is a course-relative building block of an authored account with a locally complete explanatory purpose, question, outcome, prerequisites and representations. A theorem, proof segment, example, motivation, code explanation or connection can form a unit. Syntheses likewise belong within the unit model. The definition concerns a planning resolution applicable beyond mathematics.

Initial scope is calibrated by a viable five-to-seven-minute expert exposition to the intended audience. The teacher states the assumed background and rehearses the necessary representation and justification, using notes or other prepared representations as needed. Rushing, concealed dependencies or several independent purposes indicate a design smell: revisit prerequisites or subdivide meaningfully, retaining a unit connecting the parts. This examines explanatory scope, rather than rating a lecturer's spontaneous pace. Student explanations subsequently help check these choices. Several units may compose a longer teacher-led development, including the full Sudoku account.

This heuristic has a cognitive rationale in meaningful segmentation and opportunities to process what has been presented. In probability animations, Spanjers and colleagues found that pauses improved posttest performance, while cues marking boundaries reduced reported mental effort \autocite{SpanjersEtAl2012Segmentation}. These findings motivate attention to both structure and processing opportunities; they do not determine the proposed minute range. Tabbers and de Koeijer found improved transfer under learner control alongside greater time-on-task \autocite{TabbersDeKoeijer2010Control}. A compact presentation and the activity needed to learn from it therefore warrant separate estimates.

Mathematical familiarity also changes the appropriate boundaries. Barnard and Tall's account of cognitive units and connections describes how earlier work can become one component whose internal relations remain recoverable \autocite{BarnardTall1997CognitiveUnits}. A didactic unit is a teaching specification, but its prerequisites depend on such learning. Several units may later be treated as one component, while a short expert account may initially contain several unfamiliar elements for a student. Working-memory and cognitive-load research further motivate checking these simultaneous demands and their dependence on prior knowledge \autocite{Cowan2001Capacity,ChenPaasSweller2023,DeJong2010}.

Learning-object literature emphasises reuse while considering instructional context \autocite{Wiley2000LearningObjectsOnline,Merrill2001Components}. Portability is not a defining property of a didactic unit: reuse requires reconsidering the receiving course's question, audience and dependencies. A page, developed slide or short code fragment can help inspect scope, without defining conversions between representations.

Expert calibration, student activity and total workload are distinct planning quantities. A student episode must allow exposition, dialogue and supported repair. Advance preparation, participation in other presentations and later study enter the overall workload. Allocations are revised against what students and teachers actually need to do.
\section{Three worked design analyses}
\label{sec:org7cf82c5}

The three analyses concern service mathematics for computing: logical specification, algebraic and geometric representations, and state changes. Each makes its realism lenses, ledger questions, contrasts and synthesis explicit. The Sudoku account reflects the author's teacher-led configuration; the short flawed explanations and repair prompts below are constructed illustrations of possible dialogue, not recorded classroom episodes.
\subsection{Propositional logic through compatible choices}
\label{sec:org2e1edaa}

How can a formula describe compatible choices? Scheduling and planning supply the professional purpose; variables, formulas and correctness arguments supply the disciplinary account. Experiential realism motivates retaining ordinary nine-by-nine Sudoku with three-by-three boxes as a recognisable problem. The author's teacher-led development was inspired by Raymond Hettinger's PyCon US 2019 talk and accompanying SAT notes \autocite{Hettinger2019Solvers,Hettinger2019SATNotes}. Earlier SAT encoding and computational logic teaching provide scholarly precedents \autocite{LynceOuaknine2006Sudoku,CernaEtAl2020ComputationalLogic}. The solver is supplied while the exposition opens the representation and its justification.

Let \(p,q,r\) represent three choices. Negation \(\neg\), inclusive disjunction \(\lor\) and conjunction \(\land\) express \emph{not}, \emph{or} and \emph{and}, respectively. Consider

\begin{align}
 A&=p\lor q\lor r,\\
 M&=(\neg p\lor\neg q)\land(\neg p\lor\neg r)\land(\neg q\lor\neg r),\\
 E&=A\land M.
\end{align}

The first requires at least one choice; the second forbids every simultaneous pair. Their conjunction requires exactly one, and every assignment selecting exactly one satisfies it. Cataloguing zero, one and several selections exposes what a successful example leaves undecided (Table \ref{tab:org0eb8d8d}).

\begin{table}[tbp]
\caption{\label{tab:org0eb8d8d}Acceptance of contrasting truth assignments by three restrictions; 1 means true and 0 false.}
\centering
\begin{tabular}{lccc}
Restriction & \((1,0,0)\) & \((0,0,0)\) & \((1,1,0)\)\\
\hline
\(A\): at least one & 1 & 0 & 1\\
\(M\): at most one & 1 & 1 & 0\\
\(E\): exactly one & 1 & 0 & 0\\
\end{tabular}
\end{table}

An illustrative explanation, "The disjunction says exactly one is chosen", can be tested with two true variables. Asking also about choosing nothing leads to a repaired account separating existence from uniqueness. The teacher can develop these contrasts within the exposition; subsequent student explanation can establish whether the distinction is available for use. Synthesis makes exactly one a reusable construction, while the ledger's idling question exposes a limit to its analogy with scheduling: \(M\) permits inactivity, whereas \(E\) excludes it.

A literal is a variable or its negation; a clause is a disjunction of literals. A conjunction of clauses is in conjunctive normal form (CNF), a standard SAT input form. Teaching includes deriving CNF from general formulas as well as assembling clauses directly. An example also used in Hettinger's notes is \(p\Rightarrow(q\land r)\), equivalent to \((\neg p\lor q)\land(\neg p\lor r)\) \autocite{Hettinger2019SATNotes}. The assignment \((p,q,r)=(1,1,0)\) separates this from the incorrect single clause \(\neg p\lor q\lor r\). Thus the ledger checks the restriction's meaning and the preservation of that meaning under transformation. More general cardinality constructions can extend the account \autocite{Sinz2005Cardinality}.

For full Sudoku, \(p_{rcd}\) states that cell \((r,c)\) contains digit \(d\), with all three indices ranging from 1 to 9. Exactly-one restrictions apply to each cell and each digit's occurrences in every row, column and three-by-three box; givens add single-literal clauses. A satisfying assignment determines a valid completion, and every valid completion determines a satisfying assignment. A grid satisfying row and column restrictions but violating boxes tests whether the global account is complete. Existence and uniqueness of a completion remain separate questions.

Teacher synthesis connects the local construction, its repeated use and the two directions of correctness. Analogy locates a shared pattern; deduction establishes its adequacy. Several units compose this full exposition, rather than limiting the whole development to one calibrated interval. The supplement gives a proposed eight-unit map, general cardinality formulas and verification details. Students can later reconstruct or adapt parts of the account without having to study the supplied solver's algorithm.
\subsection{Grayscale as a projection}
\label{sec:org095d40a}

Why should converting an image to grayscale leave the result unchanged on reapplication, and which differences disappear? The professional purpose concerns image processing; disciplinary realism asks how that operation is represented and explained in the relevant field; experiential realism begins with the visible change from colour to gray. The mathematical task connects those perspectives through a stated model.

Represent a pixel by normalised channels \((r,g,b)\in[0,1]^3\), with equal channels representing gray. Consider

\[
 T(r,g,b)=(s,s,s),\qquad s=0.3r+0.6g+0.1b.
\]

The model excludes rounding and nonlinear colour encoding. Its linear extension to \(\mathbb R^3\) has three identical matrix rows. Students can reconstruct them from the outputs of the standard basis vectors. Writing \(u=(1,1,1)^{\mathsf T}\) and \(w=(0.3,0.6,0.1)^{\mathsf T}\) gives \(T(x)=u(w^{\mathsf T}x)\). Since \(w^{\mathsf T}u=1\), the transformation fixes every point of the grayscale line and satisfies \(T^2=T\): it is a linear projection.

The ledger distinguishes reaching the grayscale line from fixing it, and asks which inputs become indistinguishable. An illustrative explanation, "Every output is gray, so this is a projection", fails when \(w\) is replaced by \(0.9w\): the outputs remain gray, but reapplication changes every nonzero gray input. A repair prompt asks the student to apply the map twice and identify the missing condition. The repaired account distinguishes range membership from idempotence. Separately, \((0.5,0.5,0.5)\) and \((0.7,0.4,0.5)\) have the same output; their difference lies in the kernel, the set mapped to zero.

Basis action, the kernel and idempotence offer three units. Cataloguing the contrasting maps and inputs reveals distinct explanatory obligations; abstraction and synthesis connect the coordinate recipe with its algebraic and geometric properties. The resulting transformation becomes an object for composition, with its representations and justification available for recovery. This encapsulation of mathematical work is distinct from the transformation's loss of colour distinctions.
\subsection{A contract for software failure}
\label{sec:org7c541c1}

How can an item be moved between collections without loss or duplication if an operation fails? Professional realism concerns reliable state changes; disciplinary realism concerns contracts, traces and their justification; moving an item and noticing its loss or duplication provides an experiential entry. Let the visible state consist of disjoint sets \((S,D)\), with \(x\in S\). Success must produce \((S\setminus\{x\},D\cup\{x\})\); reported failure must preserve the original pair.

An illustrative explanation, "Remove it from the source, then add it to the destination", covers success. Asking what happens if the second step fails exposes loss; reversing the order exposes duplication. Cataloguing these traces gives the ledger its contrasts: what is visible at each failure point, is the original union preserved, and are the sets disjoint at entry and on reported success or failure?

A supported repair can begin by asking which work could be completed before changing the visible state. The construction prepares the new pair in temporary copies, then installs it through a stipulated atomic replacement that either succeeds completely or leaves the state unchanged. There are no later fallible steps. The repaired explanation must cover every permitted failure and successful completion. This sequential abstract model excludes concurrency and external effects; a senior course can investigate the replacement primitive.

The contract, failure traces and repair form three possible units. Synthesis abstracts a preparation-and-commitment pattern while retaining the assumptions that make it valid. Its mathematical content consists of state relations and quantified claims about executions. The example also shows how units can organise mathematical reasoning within a computing course.
\section{Planning student and teacher effort}
\label{sec:org82fead5}

The question ledger records explanatory obligations; a unit specification connects them to prerequisites, activity and time. Table \ref{tab:org1833193} illustrates that distinction for exactly one. The unit can appear within the teacher-led Sudoku account or support a subsequent student explanation.

\begin{table}[tbp]
\caption{\label{tab:org1833193}Didactic unit specification record for the exactly-one construction.}
\centering
\begin{tabular}{p{0.23\linewidth}p{0.65\linewidth}}
Field & Specification\\
\hline
Question and outcome & Require exactly one choice and justify both directions.\\
Prerequisites & Variables, assignments, disjunction and pairwise exclusion.\\
Contrasts & No choice, one choice, two choices; required use versus allowed idling.\\
Student work & Follow and reconstruct the account; prepare an explanation when assigned, respond to questions and articulate repairs.\\
Synthesis & Separate existence from uniqueness; retain a reusable construction with stated assumptions.\\
Expert calibration & Rehearse an adequate account under the stated prerequisites.\\
Student episode & Allow explanation, dialogue and articulation of the repaired account.\\
\end{tabular}
\end{table}

Session planning reserves time for preparation, exposition, dialogue, repair and synthesis. Longer exchanges require fewer presentations or changed allocations. Synthesis can be distributed through the session, with each contribution counted once. A teacher-led development uses a different allocation of the same available time.

Student workload includes participation in the full session, advance preparation and later study. The European Credit Transfer and Accumulation System (ECTS) relates credits to outcomes and the work of learning activities \autocite{EuropeanCommission2015ECTS}. Active teaching time, scheduled commitment and independent work should be distinguished under an explicit accounting convention. Breaks and transitions occupy time too; they cannot simply be reassigned as additional study. Unit counts organise content, while activities determine workload.

Cohort size and the available sessions constrain individual presentation opportunities. Repeated attempts and teacher-led sessions must also fit the course. Larger cohorts may require parallel sessions with teaching resources and shared questions and criteria. Selected plenary contributions provide a different distribution of opportunities.

Teacher effort includes authoring, preparing examples, conducting dialogue, assessing explanations and maintaining materials. Preparation can begin with an existing account and an unresolved ledger question. An initial record identifies the intended explanation, its prerequisites, a useful contrast and the connection to retain. Once these obligations are met, further examples and extensions can be deferred. This gives preparation a starting point and a scope boundary without requiring every possible difficulty to be anticipated.

The same record can support live orientation through a compact view of the current question, essential distinctions and intended synthesis. Reuse should avoid duplicate documentation. A brief note can retain an unresolved issue for later revision; the teacher need not pursue every connection immediately. These are proposed ways of supporting sustained effort, whose usefulness must be judged against the work needed to maintain them. If preparation or follow-up becomes excessive, the design requires simplification or additional resources.
\section{Discussion}
\label{sec:orga055462}

\subsection{Design choices and tradeoffs}
\label{sec:orgb01386a}

Integration questions test whether units form a coherent account; fragmentation prompts a connecting unit or refactoring. Preparation can enable dialogue while adding work outside class. Group rehearsal uses contact time, and direct exposition establishes prerequisites and standards while reducing time for presentations. Student choice permits enrichment within stated outcomes. Planning makes these adjustments explicit.

The examples establish a computing focus without exhausting possible applications. An engineering course could organise units around a differential-equation model, contrasting a steady-state solution with a transient satisfying an initial condition. Its ledger would test assumptions, interpretation and the relation between solutions. This is a possible extension of the design procedure, requiring its own course-specific development and evaluation.
\subsection{Participation conditions for students and teachers}
\label{sec:orgf8c4a6d}

Conditions of participation belong within the design for both students and teachers. Advance materials can make expectations predictable while moving initiation, reading and time management outside class. Group preparation adds coordination and social demands; oral explanation makes spontaneous language production and public performance consequential alongside mathematical understanding.

Interviews with STEM students with ADHD and/or specific learning disabilities describe both benefits and barriers in active learning \autocite{PfeiferEtAl2023Instructor}. A survey of university students reporting autism diagnoses identifies academic strengths alongside social and learning-environment difficulties \autocite{GurbuzEtAl2019University}. Interviews with students across several neurodivergent groups also show that group work varies with its organisation and interpersonal conditions \autocite{SalvatoreEtAl2024GroupWork}. These findings motivate attention to access; they do not establish the effectiveness of the proposed method for a diagnostic group.

Adjustable conditions include advance prompts, preparation steps, individual rehearsal, notes or diagrams, pauses, predictable questions and a quieter or smaller audience. Individual preparation can remain available where group work is not an outcome. Oral exchange remains central, with assessment attending to mathematical contributions rather than rapid public participation. Decisions should respond to learners and applicable arrangements, without assigning teaching styles or mathematical directions to diagnoses.

Teacher ADHD brings the sustainability of preparation and teaching into focus. Clark and Jones report language teachers' experiences of difficulties initiating and organising work, alongside strategies such as templates and reused materials; some participants also described costly overpreparation \autocite{ClarkJones2026TeachersADHD}. Hamilton and colleagues' interviews with early-career academics show how valued variety and autonomy can coexist with overcommitment and exhaustion \autocite{HamiltonEtAl2026AcademicsADHD}. These accounts identify relevant demands and tradeoffs; they do not establish the effectiveness of this methodology for teachers with ADHD.

Applied here, these findings motivate attention to starting preparation, containing its scope and sustaining effort across teaching commitments. A bounded question can make the next preparation step explicit, while reusable units and a compact ledger can preserve earlier decisions. Prepared synthesis points allow adaptation without requiring the whole account to be composed during class. These supports must remain adjustable: a comprehensive record can itself become a maintenance burden. Evaluation should examine both the teaching achieved and its continuing cost to the teacher.
\subsection{Scope and evaluation}
\label{sec:org255f80f}

The article develops a methodology from teaching practice and examines it through worked design analyses; systematic evaluation of its educational effects remains future work. Mathematical arguments and finite checks establish properties of the examples, while learning, participation and workload require their own evaluation.

The flipped-mathematics literature gives reasons to investigate activities and conditions rather than the format alone. Maciejewski reports attainment differences in a large calculus comparison, whereas Hardebolle and colleagues' multi-year linear-algebra study finds no significant average attainment advantage and examines differing patterns by prior knowledge and gender \autocite{Maciejewski2016Flipping,HardebolleEtAl2022LinearAlgebra}. Feudel and Fehlinger provide a close precedent for expert modelling within a flipped arrangement \autocite{FeudelFehlinger2023LectureOriented}. The present contribution coordinates that function with student explanation, supported repair and course authoring.

Evaluation can time adequate expert accounts and record which prerequisites or boundaries student work exposes. It should distinguish that calibration from student activity and total workload, and examine preparation effort, reuse and maintenance alongside support, synthesis and participation. Explanations after repair, later reconstruction and adaptation address different outcomes; comparison requires attention to prior knowledge and assessment conditions.

The realism lenses and compass can be examined through the questions, contrasts and revisions they support. Compression concerns a construction's use and recoverable grounds, which brevity alone cannot establish. Ledger revisions can trace how feedback changes course design; student participation in authoring it remains a further research question. The method's central commitment is to make understanding concrete enough to plan, explain and examine as students' mathematical work becomes more connected and reusable.
\section{Acknowledgements}
\label{sec:org92b015a}

Part of this work was carried out at the Abdus Salam International Centre for Theoretical Physics (ICTP), Trieste, Italy. The work was supported by internal project funding from the Faculty of Informatics, Juraj Dobrila University of Pula.
\subsection{AI assistance}
\label{sec:org2729726}

OpenAI ChatGPT and Codex assisted literature discovery, manuscript drafting and
revision, and development of supplementary verification scripts during
November 2026. The author started with drafts and then directed the work through
written instructions and spoken comments and reviewed successive drafts;
responsibility for the content remains with the author.
\section{Ethics and availability of materials}
\label{sec:orge1996e4}

The article presents a practice-informed conceptual methodology and constructed mathematical analyses. No human-participant data were collected or analysed for this article. Illustrative student explanations and repairs are constructed examples, not transcripts or assessments of identifiable students. Descriptions of teaching arrangements do not report individual student records or establish educational effects. The supplement provides mathematical details, a proposed unit map, verification code and outputs.
\section{Competing interests}
\label{sec:org40ca5a7}

The author declares no competing interests.
\section{References}
\label{sec:orgfc3a9e7}

\printbibliography[heading=none]
\end{document}